\documentclass[11pt]{amsart}
\usepackage{amsmath,amssymb,amsfonts,url,esint}
\numberwithin{equation}{section}

\newtheorem{thm}{Theorem}[section]
\newtheorem{lem}{Lemma}[section]
\newtheorem{rem}{Remark}[section]

\begin{document}
\title[Monge-Amp\`ere equation]{ Global Monge-Amp\'ere equation with asymptotically periodic data} \subjclass{35J96; 35J67}
\keywords{Monge-Ampere equation, Dirichlet problem, a priori estimate, maximum principle, viscosity solution}
\author{Eduardo V. Teixeira}
\address{Universidade Federal do Cear\'a\\ Departamento de
Matem\'atica\\ Av. Humberto Monte, s/n, Campus do Pici - Bloco 914
\\
Fortaleza-CE, Brazil.  CEP 60.455-760} \email{teixeira@mat.ufc.br}

\author{Lei Zhang}
\address{Department of Mathematics\\
        University of Florida\\
        358 Little Hall P.O.Box 118105\\
        Gainesville FL 32611-8105}
\email{leizhang@ufl.edu}

\date{\today}

\begin{abstract} Let $u$ be a convex solution to $\det(D^2u)=f$ in $\mathbb R^n$ where $f\in C^{1,\alpha}(\mathbb R^n)$ is asymptotically close to a periodic function $f_p$. We prove that the difference between $u$ and a parabola is asymptotically close to a periodic function at infinity, for dimension $n\ge 3$.
\end{abstract}

\maketitle

\section{Introduction}

In this article we study convex, entire viscosity solutions $u \colon \mathbb R^n \to \mathbb{R}$, to the Monge-Amp\'ere equation
\begin{equation}\label{mp-g}
	\det(D^2u)=f(x),\quad \mbox{in }\quad \mathbb R^n.
\end{equation}
The forcing term $f$ is assumed to be positive and asymptotically close to a periodic function at infinity. Our main goal is to establish a classification theorem for such solutions.

\par

Monge-Amp\'ere equation with periodic data can be found in various topics in applied mathematics such as homogenization theory, optimal transportation problems, vorticity arrays, etc. Equation (\ref{mp-g}) also appears in differential geometry, when it is lifted from a Hessian manifold \cite{caf-via}. In spite of the profusion in application, Monge-Amp\'ere equation is well known for its analytical difficulty and it is no exception for equation (\ref{mp-g}) when the right hand side is close to a periodic function. In \cite{caf-li-2} Caffarelli and Li proved that if $f$ is a positive periodic function, $u$ has to be a parabola plus a periodic function with the same periodicity of $f$. This theorem can be viewed as an extension to the classification theorems of J\"orgens \cite{jorgens}, Calabi \cite{calabi}, Pogorelov \cite{p5}, Caffarelli-Li \cite{caf-li-1} for the Monge-Amp\`ere equation $det(D^2u)=1$.

The aim of this article is to establish an optimal perturbation result from the Caffarelli-Li's classification theorem, as to cover forcing terms $f$ that are asymptotically a periodic function at infinity.

In more precise terms, we let
\begin{equation}\label{aa}
 \mathcal{A} \quad \mbox{ be the set of all positive definite and symmetric matrices,}
 \end{equation}
  the assumption on $f$ is as follows:
Let $f_p$ be a positive, $C^{1,\alpha}$ periodic function in $\mathbb R^n$, i.e.:
\begin{align}\label{fp}
 &\exists d_0>0 , \alpha\in (0,1), a_1,...,a_n>0 \mbox{ such that }\\
 & d_0^{-1}\le f_p\le d_0, \quad
\|f_p\|_{C^{1,\alpha}(\mathbb R^n)}\le d_0, \nonumber\\
 & f_p(x+a_i e_i)=f_p(x), \quad \forall x\in \mathbb R^n. \nonumber
\end{align}
where $e_1=(1,0,..,0)$,..., $e_n=(0,...,0,1)$. We assume that
$f\in C^{1,\alpha}(\mathbb R^n)$ is asymptotically close to $f_p$ in the following sense:
\begin{align}\label{af}
&\exists d_1>0 \mbox{ and } \beta>2,\mbox{such that } \\
&d_1^{-1}\le f(x)\le d_1, \,\forall x\in \mathbb R^n, \quad \|f\|_{C^{1,\alpha}(\mathbb R^n)}\le d_1,
\nonumber\\
&|\nabla^j(f^{\frac 1n}-f_p^{\frac 1n})(x)|\le d_1(1+|x|)^{-\beta-j},\quad \forall x\in \mathbb R^n,\quad j=0,1,2,3.\nonumber
\end{align}
\begin{rem} $f,f_p\in C^{1,\alpha}(\mathbb R^n)$ but the difference between $f^{\frac 1n}$ and $f_p^{\frac 1n}$ is more smooth.
\end{rem}
Under the above framework, our main theorem is:
\begin{thm}\label{thm1}
Let $n\ge 3$ and $u\in C^{3,\alpha}(\mathbb R^n)$ be a convex solution to
\begin{equation}\label{mp}
\det(D^2u)=f,\quad \mbox{ in }\quad \mathbb R^n
\end{equation}
where $f$ satisfies (\ref{af}). Then there exist $b\in \mathbb R^n$, $A\in \mathcal{A}$ (defined in (\ref{aa})) with
 $det(A)=\fint_{\Pi_{1\le i\le n}[0,a_i]}f$, and $v\in C^{3,\alpha}(\mathbb R^n)$, which is $a_i-$periodic in the $i$-th variable, such that
\begin{equation}\label{mes}
	\left | u(x)- \left (\frac 12x'Ax+b\cdot x+v(x) \right ) \right |\le C \left ( 1+|x| \right )^{-\sigma},\quad \forall x\in \mathbb R^n
\end{equation}
for some $C(d_0,d_1,n,\beta, a_1,...,a_n)>0$ and $\sigma:=\min\{\beta,n-2\}$.
\end{thm}
Here we note that $\fint_{\Pi_{1\le i\le n}[0,a_i]}f$ is the average of $f$ over the box.

\begin{rem} Theorem \ref{thm1} does not imply corresponding (better) estimates on higher order derivatives because of the oscillation of $D^2u$.
\end{rem}

Caffarelli and Li \cite{caf-li-2} proved that
$$u(x)=\frac 12 x'Ax+b\cdot x+v(x), \quad \mathbb R^n,\quad n\ge 2 $$
if $f=f_p$ in (\ref{mp-g}).
Thus Theorem \ref{thm1} is an extension of the theorem of Caffarelli-Li.

The assumption $\beta>2$ is essentially optimal, as one can  observe from the following example: let $f$ be a radial, smooth, positive function such that $f(r)\equiv 1$ for $r\in [0,1]$ and $f(r)=1+r^{-2}$ for $r>2$. Let
$$u(r)=n^{\frac 1n}\int_0^r\bigg (\int_0^s t^{n-1}f(t)dt\bigg )^{\frac 1n}ds,\quad r=|x|. $$
It is easy to check that $det(D^2u)=f$ in $\mathbb R^n$. Moreover for $n\ge 3$,
$$u(x)=\frac 12|x|^2+O(\log |x|)$$
at infinity, which means by taking $f_p\equiv 1$ the estimate in Theorem \ref{thm1} is violated for $n\ge 4$.

One major difficulty in the study of (\ref{mp-g}) is that the right hand side oscillates wildly when it is scaled. The regularity theory for Monge-Amp\'ere equations with oscillating right hand side is very challenging (see \cite{figalli,savin2} ). In this respect, we found that the strategy implemented by  Caffarelli and Li in  \cite{caf-li-2} is nearly perfect, as we can only simplify a small part of their augument in \cite{caf-li-2}. In turn, our arguments are essentially based on the corresponding steps from \cite{caf-li-2} as well as previous works of Caffarelli and Li \cite{caf-ann1,caf-ann2,caf-li-1}.
 The main difference, though, is that in the proof of Theorem \ref{thm1} one needs to take care of perturbational terms in a sharp manner. In order to handle all the perturbations, we need to make use of intrinsic structures implied by Monge-Amp\'ere equations (such as
(\ref{div-free}) below), estimates on Green's functions by Littman-Stampacchia-Weinberger \cite{littman} and Krylov-Safonov Harnack inequalities, etc.

Theorem \ref{thm1} is closely related to the exterior Dirichlet problem: Given a strictly convex set $D$ and the value of $u$ on $\partial D$, can one solve the Monge-Amp\'ere equation in $\mathbb R^n\setminus D$ if the asymptotic behavior of $u$ at infinity is prescribed? Clearly Theorem \ref{thm1} must be established before such a question can be attacked. We plan to address the exterior Dirichlet problem in a future work. The traditional (interior) Dirichlet problem has been fairly well understood through the contribution of many people (see \cite{alek1, b1, n2,calabi,p5,cheng-yau1,cns,caf-ann2,ivo-1,ivo-2,ivo-3,k2,jian-wang,tw7,urbas,savin2,figalli} and the references therein). If $f$ is equal to a positive constant outside a compact set, Delano\"e \cite{delanoe}, Ferrer-Mart\'inez-Mil\'an \cite{FMM1,FMM2} and Bao-Li \cite{bao} studied the exterior Dirichlet problem for $n=2$, Caffarelli-Li \cite{caf-li-1}) studied the case of $n\ge 3$. If $f$ is a perturbation of a positive constant at infinity, Bao-Li-Zhang \cite{bao} studied the exterior Dirichlet problem in \cite{zhang-mp}.

The organization of this article is as follows: The proof of Theorem \ref{thm1} consists of five steps. First in step one we employ the argument in \cite{caf-ann2,caf-li-1} to show that the growth of the solution of (\ref{mp-g}) is roughly similar to that of a parabola. Then in step two we prove that $D^2u$ is positive definite, which makes (\ref{mp-g}) uniformly elliptic. The key point in this step is to consider a second order incremental of $u$ as a subsolution to an elliptic equation. In step three we prove a pointwise estimate of the second incremental of $u$. The proof of Theorem \ref{thm1} for $n\ge 4$ is placed in step four since all the perturbations in this case are bounded. Finally in step five we prove the case $n=3$, which is a little different because of a logarithmic term. We use the Krylov-Safonov Harnack inequality for linear equations to overcome the difficulties caused by the logarithmic term.

\medskip

\noindent {\bf Acknowledgement:} The authors would like to thank the hospitality of Universidade Federal do Cear\'a and University of Florida, where this work was conducted. ET acknowledges support from CNPq and Funcap. Both authors are grateful to the anonymous referee, who read the whole paper very carefully and gave many excellent suggestions.

\section{Proof of Theorem \ref{thm1}}

Since the Monge-Amp\'ere equation is invariant under affine transformation, we assume $a_1=...=a_n=1$ and $\fint_{[0,1]^n} f=1$.

In step one we prove that $u$ grows like a quadratic polynomial at infinity. First we normalize $u$ to make $u(0)=0$ and $u\ge 0$ in $\mathbb R^n$. Since $f$ is bounded above and below by two positive constants, we use the argument in Caffarelli-Li \cite{caf-li-1}, see also \cite{zhang-mp}. Let
$$\Omega_M=\{x\in \mathbb R^n, \quad u(x)<M\}. $$

Then the following properties hold:
\begin{enumerate}
\item $C^{-1}M^{\frac n2}\le |\Omega_M |\le CM^{\frac n2}$,
\item There exists $A_M(x)=a_M\cdot x+b$ such that $det(a_M)=1$ and $B_R\subset A_M(\Omega_M)\subset B_{nR}$ and
\item $\frac 1C\sqrt{M}\le R\le C\sqrt{M}$,
\item $2nR\ge dist(a_M(\Omega_{\frac M2}),\partial a_M(\Omega_M))\ge \frac RC$
\end{enumerate}
where all the constants $C$ only depend on $d_1$ and $n$.

All the properties listed above are proved in \cite{caf-li-1} only based on the assumption that $f$ is bounded above and below by two positive constants. Here for the convenience of the reader we mention the idea of the proof:
First Caffarelli-Li used the following lemma (Lemma 2.9 in \cite{caf-li-1})

\medskip

\emph{Lemma A (Caffarelli-Li):  Let $e_1=(1,0,..,0)$ and
$$B_{\delta}'=\{(0,x_2,...,x_n);\quad |(0,x_2,...,x_n)|<\delta \}. $$
 Let $K$ be the convex hull of $\bar B_{\delta}'\cup \{re_1\}$, and let $u$ be a nonnegative convex viscosity solutions of
$det(D^2u)\ge \lambda>0$ in the interior of $K$.
Assume
$u\le \beta$ on $\bar B_{\delta}'. $
Then there exists $C(n)>1$ such that
$$\max \{\beta, u(re_1)\} \ge \frac{\lambda^{1/n}\delta^{\frac{2(n-1)}n}r^{2/n}}{C}. $$
}

In other words, $u$ cannot be small in one direction for too long. Lemma A is important since it says Pogorelov's famous example of non-strictly convex solution does not exist if the domain is large. A simple application of Lemma A leads to $\Omega_M\subset B_{CM^{\frac n2}}$. A volume preserving affine transformation can be used to make the image of $\partial \Omega_M$ between too balls with comparable radii. A comparison with a parabola gives $|\Omega_M|\sim M^{\frac n2}$ using only the upper bound and lower bound of $f$.  In \cite{caf-ann1} Caffarelli proved that $u$ must depart from its level set in a non-tangential manner, using this we have
$dist(A_M(\partial \Omega_M),A_M(\Omega_{M/2}))\sim M^{\frac 12}$ where $A_M$ is a volume-preserving affine transformation:
$$A_M(x)=a_Mx+b_M,\quad det(a_M)=1. $$
Using $u(0)=0$ and $u\ge 0$ one can further conclude
\begin{equation}\label{like-ball}
B_{R/C}\subset a_M(\Omega_M)\subset B_{2nR}\quad \mbox{ where } \quad R=M^{\frac 12}.
\end{equation}
Equation (\ref{like-ball}) makes it convenient to define
$$u_M(x)=\frac 1{R^2}u(a_M^{-1}(Rx)),\quad x\in O_M:=\frac 1R a_M(\Omega_M), $$
Using $u(0)=0$ and $u\ge 0$ we have
$$B_{1/C}\subset O_M\subset B_{2n}. $$
In order to analyze the level surfaces of $u_M$, which satisfies
\begin{equation}\label{mp-m}
det(D^2u_M(\cdot))=f(a_M^{-1}(R\cdot)),\quad \mbox{ in }\quad O_M.
\end{equation}
We use the following equation to give a good approximation of $u_M$:
\begin{equation}\label{mp-m2}
\left\{\begin{array}{ll}
det(D^2w_p)=f_p(a_M^{-1}(Rx)),\quad \mbox{ in } O_M, \\
w_p=M/R^2,\quad \mbox{ on }\quad \partial O_M.
\end{array}
\right.
\end{equation}
Let $h=u_M-w_p$ be the difference of $u_M$ and $w_p$, then by the Alexandrov estimate ( see \cite{caf-li-1,zhang-mp}) we have
$$\max_{O_M}(h^-)\le C\bigg (\int_{S^+}det(D^2(u_M-w_p))\bigg )^{\frac 1n}$$
where $h^-$ is the negative part of $h$: ($h=h^+-h^-$) and
$$S^+=\{ x\in O_M;\quad D^2(u_M-w_p)>0\quad \}. $$
On $S^+$ by the concavity of $det^{\frac 1n}$ on positive definite matrices we have
$$det^{\frac 1n}(\frac{D^2u_M}2)\ge \frac 12det^{\frac 1n}(D^2(u_M-w_p))+\frac 12det^{\frac 1n}(D^2w_p). $$
Thus
$$\max_{O_M}h^-\le C\bigg (\int_{O_M}|f^{\frac 1n}(a_M^{-1}(Rx))-f_p^{\frac 1n}(a_M^{-1}(Rx))|^ndx\bigg )^{\frac 1n}. $$
It is easy to see from (\ref{af}) that the right hand side is $O(1/R)$. Next we cite the homogenization theorem of Caffarelli-Li (Theorem 3 of \cite{caf-li-2}): Let $w$ satisfy
$$det(D^2w)=1,\quad \mbox{ in } \quad O_M$$ with $w=M/R^2$ on $\partial O_M$, then
$$|w-w_p|\le CR^{-\delta} $$
for some $\delta (d_0,d_1)>0$. Thus
\begin{equation}\label{close-1}
|u_M(x)-w(x)|\le CR^{-\delta}
\end{equation}
for some $\delta>0$. Set
$$E_M:=\{x;\quad (x-\bar x)'D^2w(\bar x)(x-\bar x)\le 1\} $$
where $\bar x$ is the unique minimum of $w$ (note that Caffarelli \cite{caf-ann1} proved that the minimum point of $w$ is unique) that satisfies $dist(\bar x,\partial O_M)>C(d_1,n)$. By the same argument in \cite{caf-ann2,caf-li-1} there exist $\bar k$ and $C$ depending only on $n$ and $d_1$ such that for $\epsilon<<\delta$, $M=2^{(1+\epsilon)k}$, $2^{k-1}\le M'\le 2^k$, $R\sim M^{\frac 12}$,
$$(\frac{2M'}{R^2}-C2^{-\frac{3\epsilon k}2})^{\frac 12}E_M\subset
\frac 1Ra_M(\Omega_{M'})\subset (\frac{2M'}{R^2}+C2^{-\frac{3\epsilon k}2})^{\frac 12}E_M,\quad \forall k\ge \bar k, $$
which can be translated as
$$\sqrt{2M'}(1-\frac{C}{2^{\epsilon k/2}})E_M\subset a_M(\Omega_{M'})\subset \sqrt{2M'}(1+\frac{C}{2^{\epsilon k/2}})E_M. $$
Let $Q$ be a positive definite matrix satisfying $Q^2=D^2w(\bar x)$, $O$ be an orthogonal matrix that makes $T_k=OQa_M$ upper triangular, then $det(T_k)=1$ and by Proposition 3.4 of \cite{caf-li-1}
$$\|T_k-T\|\le C2^{-\frac{\epsilon k}2} $$
for some matrix $T$. By setting
$$v=u\cdot T$$
we have
$$det(D^2 v(x))=f(Tx) $$
and
$$\sqrt{2M'}(1-\frac{C}{2^{\epsilon k/2}})B_1\subset \{x;\quad v(x)<M'\}\subset \sqrt{2M'}(1+\frac{C}{2^{\epsilon k/2}})B_1 $$
for all $M'\ge 2^{\bar k}$. Here $B_1$ stands for the unit ball (throughout the article we use $B(p,r)$ to denote the ball centered at $p$ with radius $r$. If $p$ is the origin we may use $B_r$). Consequently
$$|v(x)-\frac 12 |x|^2|\le C|x|^{2-\epsilon}. $$
The equation for $v$ is
$$det(D^2 v)=f_v(x) $$
where $f_v(x)=f(T(x))$, correspondingly we let $f_{v,p}(x)=f_p(T(x))$.

\medskip
\noindent{\bf Step two: Uniform Ellipticity}
\medskip

The purpose of this step is to show: There exist $c_1$ and $c_2$ depending only on $d_0,d_1,\beta$, $a_1,...,a_n$ and $n$ such that
\begin{equation}\label{uni-ell}
c_1 I\le D^2 v\le c_2 I.
\end{equation}

First we choose $M>100$ so that for $R=M^{\frac 12}$ and
$$v_R(y)=\frac 1{R^2}v(Ry), $$
$$\Omega_{1,v_R}:=\{y;\quad v_R(y)\le 1\} $$
is very close to $B_{\sqrt{2}}$ in the sense that $\partial \Omega_{1,v_R}\subset B_{\sqrt{2}+\epsilon}-B_{\sqrt{2}-\epsilon}$ for some $\epsilon>0$ small. Applying the standard interior estimate for Monge-Amp\'ere equations (\cite{caf-ann2,jian-wang} ) we have
$$\|D^2v_R\|_{L^{\infty}(\Omega_{1,v_R})}=\|D^2v\|_{L^{\infty}(B_{2})}\le C. $$
In general for $|x|>100$ we shall prove (\ref{uni-ell}) for $D^2v(x)$.
To this end we consider the following second order incremental for $v$:
$$\Delta_e^2v(x):=\frac{v(x+e)+v(x-e)-2v(x)}{\|e\|^2} $$
where $e\in \mathbb R^n$ and $\|e\|$ is its Euclidean norm.
Later we shall always choose $e\in E$ which is defined as
\begin{align}\label{E-def}
E:=&\{a_1 v_1+...+a_nv_n;\quad a_1,...a_n\in \mathbb Z,\\
&\quad f_{v,p}(x+v_i)=f_{v,p}(x), \, \forall x\in \mathbb R^n,\quad i=1,...,n. \} \nonumber
\end{align}

Let $R_1=|x|$ (recall $|x|>10$), $e_x=\frac{x}{R_1}$ and
$$v_{R_1}(y)=\frac 1{R_1^2}v(R_1y),\quad y\in B(e_x,\frac 23). $$
 Then by the closeness result of step one, the sections: $S_{v_{R_1}}(e_x,\frac 14)$ and
$S_{v_{R_1}}(e_x,\frac 18)$ around $e_x$ are very close to the corresponding sections of the parabola $\frac 12 |y|^2$. Here we recall that for a convex, $C^1$ function $v$,
$$S_v(x, h):=\{y;\quad v(y)\le v(x)+\nabla v(x)\cdot (y-x)+h,\,\, \}. $$
In other words, $S_{v_{R_1}}(e_x,\frac 14)$ looks very similar to $B(e_x,1/\sqrt{2})$ and $S_{v_{R_1}}(e_x,\frac 18)$ is very close to $B(e_x,1/2)$.
Moreover the equation for $v_{R_1}$ is
\begin{equation}\label{vr2}
det(D^2v_{R_1}(y))=f_v(R_1y)
\end{equation}
where $f_v(R_1y)=f(T(R_1y))$.
Let $e_{R_1}=e/R_1$, then direct computation shows
\begin{align*}
&\Delta_e^2v(x)=\frac{v(x+e)+v(x-e)-2v(x)}{\|e\|^2},\quad x=R_1y\\
=&\frac{(v(R_1y+e)+v(R_1y-e)-2v(R_1y))}{\|e\|^2}\\
=&\frac{R_1^2(v_{R_1}(y+e_{R_1})+v_{R_1}(y-e_{R_1})-2v_{R_1}(y))}{\|e\|}=\Delta_{e_{R_1}}^2v_{R_1}(y).
\end{align*}
Let
$$w(y)=\frac{v_{R_1}(y+e_{R_1})+v_{R_1}(y-e_{R_1})}2 $$
and $F=det^{\frac 1n}$. Then the concavity of $F$ on positive definite matrices gives
\begin{align}
det^{\frac 1n}(D^2w)(y)&\ge \frac 12 det^{\frac 1n}(D^2 v_{R_1}(y+e_{R_1}))+\frac 12 det^{\frac 1n}(D^2v_{R_1}(y-e_{R_1}))\nonumber\\
&=\frac 12 f_v^{\frac 1n}(R_1y+e)+\frac 12 f_v^{\frac 1n}(R_1y-e). \label{det-con-1}
\end{align}
On the other hand
\begin{align} \label{det-con-2}
F(D^2w)&\le F(D^2 v_{R_1})+F_{ij}(D^2 v_{R_1})\partial_{ij}(w-v_{R_1})\\
&=f_v^{\frac 1n}(R_1y)+F_{ij}(D^2v_{R_1})\partial_{ij}(w-v_{R_1}) \nonumber
\end{align}
where
$$F_{ij}(D^2v_{R_1})=\frac{\partial F(D^2v_{R_1})}{\partial_{ij}v_{R_1}}=\frac 1n (det(D^2v_{R_1}))^{\frac 1n-1}cof_{ij}(D^2v_{R_1}).$$
Thus the combination of (\ref{det-con-1}) and (\ref{det-con-2}) gives
$$a_{ij}\partial_{ij}(\Delta_{e_{R_1}}^2 v_{R_1})\ge E_1$$
where
\begin{align*}
a_{ij}&=cof_{ij}(D^2v_{R_1}), \\
E_1:&=ndet(D^2v_{R_1})^{(n-1)/n}\\
&\cdot (\frac{R_1^2}{2\|e\|^2}(f_v^{\frac 1n}(R_1y+e)+f_v^{\frac 1n}(R_1y-e)-2f_v^{\frac 1n}(R_1y)))
\end{align*}
For applications later we state the following fact: For any smooth $u$
\begin{equation}\label{div-free}
\partial_i(cof(D^2u)_{ij})=0.
\end{equation}

Let $f_{v,p}(y)=f_p(Ty)$, then by choosing $e\in E$ (see (\ref{E-def}))
\begin{equation}\label{per-e}
f_{v,p}^{\frac 1n}(R_1y+e)+f_{v,p}^{\frac 1n}(R_1y-e)-2f_{v,p}^{\frac 1n}(R_1y)=0,\quad \mbox{ for all }y\in \mathbb R^n.
\end{equation}

By (\ref{af}) we see that, for $y\in B(e_x, \frac 14)$ (which implies $2>|y|>1-1/\sqrt{2}$) and $e\in E$ with $\|e\|\le \frac 1{10}R_1$
\begin{align}\label{error-1}
E_1(y)&=n(det(D^2v_{R_1})^{(n-1)/n})\frac{R_1^2}{2\|e\|^2}\{(f_v^{\frac 1n}(R_1y+e)-f_{v,p}^{\frac 1n}(R_1y+e))\nonumber\\
&+(f_v^{\frac 1n}(R_1y-e)-f_{v,p}^{\frac 1n}(R_1y-e))-2(f_v^{\frac 1n}(R_1y)-f_{v,p}^{\frac 1n}(R_1y))\}\nonumber\\
=&O(R_1^{-\beta}).\nonumber
\end{align}

Now we construct a function $h_0$ that solves
$$\left\{\begin{array}{ll}
a_{ij}\partial_{ij}h_0=-E_1,  \quad \mbox{ in  } B(e_x,\frac 14) \\
h_0=0\quad \mbox{ on }\quad \partial B(e_x, \frac 14)
\end{array}
\right.
$$

Then we use the following classical estimate of Aleksandrov (see Page 220-222 of \cite{gilbarg} for a proof):

\emph{Theorem A: Let $\Omega$ be a domain in $\mathbb R^n$ and let $v$ be a solution in $\Omega$ of the equation
$$a^*_{ij}\partial_{ij}v=g $$
such that $v=0$ on $\partial \Omega$ and the coefficient matrix $(a^*_{ij})_{n\times n}$ satisfies
$$c_1\le det(a^*_{ij})\le c_2, \quad \mbox{and}\quad (a^*_{ij})>0, $$
 then
$$|v(x)|\le C(n,c_1,c_2) diam(\Omega) \|g\|_{L^n(\Omega)}, \quad \forall x\in \Omega$$ }

Applying Theorem A to $h_0$ we have
\begin{equation}\label{error-2}
|h_0(y)|\le C(n,d_0)R_1^{-\beta},\quad \mbox{ for }\quad y\in B(e_x,\frac 14).
\end{equation}
\begin{rem} The estimate of Theorem A does not depend on constants of uniform ellipticity.
\end{rem}

Thus $\Delta_{e_{R_1}}^2 v_{R_1}+h_0$ is a super solution:
$$a_{ij}\partial_{ij}(\Delta_{e_{R_1}}^2 v_{R_1}+h_0)\ge 0,\quad \mbox{in}\quad B(e_x,\frac 14). $$
In order to obtain a pointwise estimate for $\Delta_{e_{R_1}}^2 v_{R_1}$ we recall the following Harnack inequality of Caffarelli-Gutierrez in \cite{gute}:

 \medskip

 \emph{Theorem B (Caffarelli-Gutierrez): Let $O$ be a convex set in $\mathbb R^n$ and $B_1\subset O\subset B_n$ ($n\ge 2$). Suppose $\phi\in C^2(O)$ satisfies, for $0<\lambda<\Lambda<\infty$,
 $$ \lambda<det(D^2\phi)\le \Lambda, \mbox{ in } \,\, O, \mbox{ and } \quad \phi=0 \mbox{ on } \quad \partial O. $$
 Then for any $r>s>0$ and any $w\in C^2(O)$ satisfying
 $$a_{ij}\partial_{ij}w\ge 0,\quad w\ge 0,$$
 where $a_{ij}=det(D^2\phi)\phi^{ij}$ and $(\phi^{ij})_{n\times n}=(D^2\phi)^{-1}$, we have
 $$\max_{x\in O, dist(x,\partial O)>r}w\le C\int_{x\in O,dist(x,\partial O)>s}w, $$
 for some $C>0$ depending only on $\lambda, \Lambda, r,s$. }

 \medskip

Using Theorem B for $\Delta_{e_{R_1}}^2 v_{R_1}+h_0$ we have
$$\max_{y\in S(v_{R_1},e_x, \frac 1{16})}(\Delta_{e_{R_1}}^2 v_{R_1}+h_0)\le C(n,d_0)\int_{S(v_{R_1},e_x,\frac 18)}(\Delta_{e_{R_1}}^2
v_{R_1}+h_0) $$
Note that the distance between the two sections above is comparable to $1$. It is also important to point out that the estimates in Theorem B does not depend on the regularity of $v_{R_1}$.
 Before we proceed we cite the following Calculus lemma that can be found in \cite{caf-li-2}:

\medskip

\emph{Lemma B: Let $g\in C^2(-1,1)$ be a strictly convex function, and let $0<|h|\le \epsilon$. Then
$$\Delta_h^2g(x)>0,\, \forall |x|\le 1-2\epsilon, \quad \mbox{ and } \quad \int_{-1+2\epsilon}^{1-2\epsilon}\Delta_h^2g\le \frac{C}{\epsilon} osc_{(-1,1)}g, $$
where $C$ is a universal constant and $osc$ stands for oscillation. }

\medskip

 Let $L$ be a line parallel to $e$, then by Lemma B,
$$\int_{S(v_{R_1},e_x,\frac 18)\cap L}\Delta_{e_{R_1}}^2 v_{R_1}\le C. $$
Consequently, letting $L$ go through all directions, we have
$$\int_{S(v_{R_1},e_x,\frac 18)}\Delta_{e_{R_1}}^2 v_{R_1}\le C. $$
Since $h_0=O(R_1^{-\beta})$  we have proved
\begin{equation}\label{quo-2}
0\le \Delta_e^2v(x)=\Delta_{e_{R_1}}^2 v_{R_1}(y)\le C,\quad \mbox{ for } y\in B(e_x,\frac 1{16}),\,\, x=R_1y.
\end{equation}
Note that $\Delta_e^2v\ge 0$ because $v$ is convex. Then the same argument as in \cite{caf-li-2} can be employed to prove that the level surfaces of $v$ are like balls.

For the convenience of the readers we describe the outline of this argument. Given $x\in \mathbb R^n$ with $|x|>100$ we set
$$\gamma:=\sup_{e\in E, \|e\|\le \frac 1{10}\|x\|}\sup_{y\in \mathbb R^n}\Delta_e^2v(y)$$
 and
$$\bar v(z)=v(z+x)-v(x)-\nabla v(x)z.$$
Clearly $\bar v(0)=0=\min_{\mathbb R^n} \bar v$.
By (\ref{quo-2}) it is easy to see
$$\sup_{B_r}\bar v\le C(n)\gamma r^2,\quad \forall r \in (1, \frac{1}{10}\|x\|).  $$
On the other hand for any $\bar z\in \partial B_r$ we show that for $r$ large (but still only depending on $n$ and $d_0$), $\bar v(\bar z)\ge 1$. Indeed, by (\ref{quo-2}),
\begin{equation}\label{uni-2}
\bar v(\frac{\bar z}2+e)+\bar v(\frac{\bar z}2-e)-2\bar v(\frac{\bar z}2)\le \gamma \|e\|^2
\end{equation}
for all $e\in E$ and $\|e\|\le \frac 1{10} \|x\|$. Hence for $z\in \frac{\bar z}2+(-2,2)^n$, we find $e\in E$ with $\|e\|\le C(n)$ so that $z$ is on the same line with $\frac{\bar z}2+e$ and $\frac{\bar z}2-e$
and is between them. Thus(\ref{uni-2}) implies
$$\bar v(z)\le \bar v(\frac{\bar z}2+e)+\bar v(\frac{\bar z}2-e)\le 2\bar v(\frac{\bar z}2)+C(n)\gamma. $$
Further more, by $\bar v(0)=0$ and the convexity of $\bar v$ we have
$$2\bar v(\frac{\bar z}2)\le \bar v(\bar z). $$
Therefore the following holds:
$$\bar v(z)\le \bar v(\bar z)+C(n)\gamma, \quad z\in \frac{\bar z}2+(-2,2)^n. $$
Consider $$w(z)=\frac{\bar v(\frac{\bar z}2+z)}{\bar v(\bar z)+C(n)\gamma}. $$
Clearly $w$ satisfies $det(D^2w)\ge d_0^{-1}/(\bar v(\bar z)+C(n)\gamma)^n$ and
$$w(z)\le 1, \quad \mbox{ for  }\quad z\in \bar z/2+(-2,2)^n. $$
Applying Lemma A  we have
$$\max \{1,\bar v(\bar z)\}\ge \bigg (\frac{d_0^{-\frac 1n} C(n)}{\bar v(\bar z)+C(n)\gamma} \bigg ) r^{\frac 2n}.  $$
Here we used the fact that $\bar v(0)=0$ and $\bar v(\bar z)$ is the maximum value of $\bar v$ on the line segment connecting $0$ and $\bar z$.
If $\bar v(\bar z)\ge \gamma$, no need to do anything ( here we assume $\gamma>1$ without loss of generality), otherwise we obviously have
$$\max\{1,\bar v(\bar z)\}\ge \frac{d_0^{-\frac 1n}}{C(n)\gamma }r^{\frac 2n}. $$
Choose $r$ large but only depending on $n$,$\gamma$ and $d_0$ we can still make the right hand side of the above greater than $1$. Therefore we have proved $\bar v(z)\ge 1$.
By standard argument the sections of $\bar v$:
$S(\bar v,0,h)$ are similar to balls if $h\sim 1$ and $h<1$.
Then interior estimate of Caffarelli, Jian-Wang \cite{jian-wang} can be employed on $S(\bar v, 0, h)$ to obtain the $C^{2,\alpha}$ norm of $\bar v$. In particular
$|D^2v(x)|=|D^2\bar v(x)|\le C$. Once the upper bound is obtained we also have the lower bounded from the Monge-Amp\'ere equation. (\ref{uni-ell}) is established.

\begin{rem} In order to prove (\ref{uni-ell}) the assumptions on the derivatives of $f^{\frac 1n}-f_p^{\frac 1n}$ are not essential. In fact, as long as
$$|(f^{\frac 1n}-f_p^{\frac 1n})(x)|\le C(1+|x|)^{-\beta}, \quad \mbox{ for }\quad x\in \mathbb R^n, $$
(\ref{uni-ell}) still holds.
\end{rem}

\noindent{\bf Step three. Pointwise estimate of $\Delta_e^2v$}

Let $e\in E$, recall the equation for $\Delta_e^2v$ is
\begin{equation}\label{p3e1}
\tilde a_{ij}(x)\partial_{ij}\Delta_e^2v\ge E_2.
\end{equation}
where
$$ \tilde a_{ij}(x)=cof_{ij}(D^2v(x))  $$ is uniformly elliptic and divergence free,
$$E_2=n(det(D^2v)^{\frac{n-1}n})\frac{f_v^{\frac 1n}(x+e)+f_v^{\frac 1n}(x-e)-2f_v^{\frac 1n}(x)}{\|e\|^2}. $$
By the assumption on $f$,
\begin{align*}
|E_2(x)|&=n(det(D^2v)^{\frac{n-1}n})\\
&\bigg(\frac{(f_v^{\frac 1n} (x+e)-f_{p,v}^{\frac 1n}(x+e))+(f_v^{\frac 1n}(x-e)-f_{v,p}^{\frac 1n}(x-e))}{\|e\|^2}\nonumber\\
&-2\frac{(f_v^{\frac 1n}(x)-f_{v,p}^{\frac 1n}(x))}{\|e\|^2}\bigg ).
\end{align*}
Here we have the following important observation:
\begin{equation}\label{e2e}
|E_2(x)|\le C(1+|x|)^{-2-\beta},\quad \forall e\in E
\end{equation}
where $C$ is independent of $e\in E$ and $x$.
Indeed, let $g(x)=f_v^{\frac 1n}(x)-f_{v,p}^{\frac 1n}(x)$, then
$$|E_2(x)|=(g(x+p)+g(x-p)-2g(x))/\|e\|^2. $$
Without loss of generality we assume that $e=(\|e\|,0,...,0)$. Then
\begin{align*}
&E_2(x)=\frac{(g(x+e)-g(x))-(g(x)-g(x-e))}{\|e\|^2}\\
=&\frac{\int_0^{\|e\|}\partial_1g(x_1+t,x')dt-\int_0^{\|e\|}\partial_1g(x_1-\|e\|+t,x')dt}{\|e\|^2}\\
=&\frac{\int_0^{\|e\|}\int_0^{\|e\|}\partial_{11}g(x_1-\|e\|+t+s,x')dtds}{\|e\|^2}\\
=&\int_0^1\int_0^1\partial_{11}g(x_1-\|e\|+\|e\|t+\|e\|s,x')dtds.
\end{align*}
Using
$$|D^2g(x)|\le C(1+|x|)^{-\beta-2},$$
which comes from (\ref{af}) we have
$$|E_2(x)|\le C\int_0^1\int_0^1\bigg ( \big (x_1+\|e\|(-1+t+s)\big )^2+|x'|^2\bigg )^{-\frac{\beta+2}2}dtds. $$
If $|x'|> \frac 18 |x|$ , (\ref{e2e}) holds obviously. So we only consider when $|x_1|> \frac 12 |x|$.
\begin{align*}
&\int_0^1\int_0^1\bigg (\big (x_1+\|e\|(-1+t+s) \big )^2+|x'|^2\bigg )^{-\frac{\beta+2}2}dtds\\
\le & \int_0^1\int_0^1 \bigg | x_1+\|e\|(-1+t+s) \bigg |^{-\beta -2}dtds \\
=&\|e\|^{-\beta -2}\int_0^1\int_0^1 |L-1+t+s |^{-\beta-2}dtds,\quad L=x_1/\|e\|.
\end{align*}
If $|L|>8$, $|L-1+t+s|>|L|/2$, which gives
\begin{align*}
&\int_0^1\int_0^1\bigg (\big (x_1+\|e\|(-1+t+s) \big )^2+|x'|^2\bigg )^{-\frac{\beta+2}2}dtds\\
 &\le C\|e\|^{-\beta -2}|L|^{-\beta-2}\le C\|x_1\|^{-\beta}\le C|x|^{-\beta-2}.
\end{align*}
where all the constants are absolute constants. Clearly (\ref{e2e}) holds in this case.

If $|L|<8$, easy to see that
$$\int_0^1\int_0^1|L-1+t+s|^{-\beta-2}dtds \le C(\beta). $$
Then
\begin{align*}
&\int_0^1\int_0^1\bigg (\big (x_1+\|e\|(-1+t+s) \big )^2+|x'|^2\bigg )^{-\frac{\beta+2}2}dtds\\
\le &C(\beta)\|e\|^{-\beta-2}\le C(\beta)\|x\|^{-\beta-2},
\end{align*}
since $|L|<8$ implies $\|e\|>|x_1|/8>|x|/16$.
 (\ref{e2e}) is established.

By standard elliptic estimates $f\in C^{1,\alpha}$ yields $\tilde a_{ij}\in C^{1,\alpha}$. Moreover $\tilde a_{ij}$ is divergence free. The Green's function $G(x,y)$ corresponding to
$-\partial_{x_i}(\tilde a_{ij}\partial_{x_j})$ satisfies
$$-\partial_{x_i}(\tilde a_{ij}\partial_{x_j}G(x,y))=\delta_y $$
and
\begin{equation}\label{stand-g}
0\le G(x,y)\le C|y-x|^{2-n}, \quad x,y\in \mathbb R^n
\end{equation}
by a result of
 Littman-Stampacchia-Weinberger in \cite{littman}.

Let $E_2^-\ge 0$ be the negative part of $E_2$: $E_2=E_2^+-E_2^-$ and
$$h(x)=-\int_{\mathbb R^n}G(x,y)E_2^-(y)dy.  $$
Then $h(x)$ satisfies
\begin{equation}\label{p3e2}
\partial_{x_i}(\tilde a_{ij}\partial_{x_j})h=E_2^-,\quad \mbox{in}\quad \mathbb R^n.
\end{equation}
Here we claim that $h$ satisfies
\begin{equation}\label{p3h1}
0\le -h(x)\le C(1+|x|)^{-\beta_1}, \quad \beta_1=\min\{n-2,\beta\}
\end{equation}

The estimate of (\ref{p3h1}) is rather standard, for $x\in \mathbb R^n$, we divide $\mathbb R^n$ into three regions:
\begin{align*}
\Omega_1=\{y;\,\, |y-x|<|x|/2\,\}\\
\Omega_2=\{y;\quad |y|<|x|/2\,\, \}, \\
\Omega_3=\mathbb R^n\setminus (\Omega_1\cup \Omega_2).
\end{align*}

Then one observes $\int_{\Omega_i}G(x,y)E_2^-(y)dy$ for each $\Omega_i$ can be estimated easily. (\ref{p3h1}) is established.
From the definition of $h$ we see that
\begin{equation}\label{delta-vij}
\tilde a_{ij}\partial_{ij}(\Delta_e^2v+h)\ge 0,\quad \mbox{ in }\quad \mathbb R^n.
\end{equation}
The main result of this step is:
\begin{equation}\label{2-e}
\mbox{ Given  }e\in E;\quad \Delta_e^2v(x)\le 1-h(x),\quad \mbox{ for }\quad x\in \mathbb R^n.
\end{equation}

Let $v_{\lambda}(x)=v(\lambda x)/\lambda^2$ and $P(x)=\frac 12 |x|^2$. First we claim that
\begin{equation}\label{2-e1}
D^j(v_{\lambda}-P(x))\to 0,\quad j=0,1,\quad \forall x\in K\subset\subset \mathbb R^n, \mbox{  as  }\lambda\to \infty,
\end{equation}
where $K$ is any fixed compact subset of $\mathbb R^n$.

To see (\ref{2-e1}), by step one
\begin{equation}\label{2-e2}
|v_{\lambda}(x)-P(x)|\le C\lambda^{-\epsilon},\quad \forall x\in K\subset\subset \mathbb R^n.
\end{equation}
and $|D^2v_{\lambda}(x)|\le C$ for $x\in K$ we obtain by Ascoli's theorem that $\partial_lv_{\lambda}(x)$ ($l=1,..,n$) tends to a continuous function. By (\ref{2-e2}) this function has to be $x_l$. Thus (\ref{2-e1}) holds.
\begin{rem} We don't have the estimate of $\|Dv_{\lambda}-x\|_{L^{\infty}(K)}$ but we don't need it.
\end{rem}

The proof of (\ref{2-e}) is as follows.

 Let $\alpha=\sup_{\mathbb R^n}(\Delta^2_ev+h)$ for $e\in E$. By step two $\alpha<\infty$. Let $\hat e=e/\lambda$, then
by (\ref{2-e1}) and  Lebesgue's dominated convergence theorem
\begin{equation}\label{caf-li-e}
\lim_{\lambda\to \infty}\int_{B_1}\Delta_{\hat e}^2v_{\lambda}=\int_{B_1} 1dx=|B_1|.
\end{equation}
Indeed, the integral over $B_1$ can be considered as the collection of integration on segments all in the direction of $\hat e$. Since $Dv_{\lambda}\to DP$ in $C^0$ and $DP$ is smooth,
the following Lemma C in \cite{caf-li-2} implies (\ref{caf-li-e}):

\medskip

\emph{Lemma C (Caffarelli-Li) Let $g_i$ converge to $g$ in $C^1[-1,1]$, $g\in C^2(-1,1)$ and $|h_i|\to 0$ as $i\to \infty$. Then for all $-1<a<b<1$,
$$\lim_{i\to \infty} \int_a^b \Delta_{h_i}^2g_i=g'(b)-g'(a)=\int_a^b g''. $$
}

\medskip

Let
$$h_{\lambda}(x)=h(\lambda x)/\lambda^2. $$
It follows immediately from (\ref{p3h1}) that
$$\lim_{\lambda\to \infty}\int_{B_1}(\Delta_{\hat e}^2v_{\lambda}+h_{\lambda})=|B_1|, $$
which obviously implies $\alpha\ge 1$ because otherwise it is easy to see that the left hand side of the above is less than $|B_1|$.  Our goal is to show that $\alpha=1$. If this is not the case we have $\alpha>1$. Then
$$\lim_{\lambda}\sup(\frac{\alpha+1}2|\{\Delta_{\hat e}^2v_{\lambda}+h_{\lambda}\ge \frac{\alpha+1}2\}\cap B_1|\le \lim_{\lambda\to \infty}\int_{B_1}\Delta_{\hat e}^2v_{\lambda}+h_{\lambda}=|B_1|. $$

Thus
$$\frac{|\{\Delta_{\hat e}^2 v_{\lambda}+h_{\lambda}\ge \frac{\alpha+1}2\}\cap B_1|}{|B_1|}\le 1-\mu $$
for $\mu=\frac 12(1-\frac{2}{\alpha+1})>0$. Here we emphasize that it is important to have $\mu>0$. Equivalently for large $\lambda$
$$\frac{|\{\Delta_{\hat e}^2v_{\lambda}+h_{\lambda}\le \frac{\alpha+1}2\}\cap B_1|}{|B_1|}\ge \mu. $$

Then we cite the following well known result which can be found in \cite{caf-cab} (Lemma 6.5):

\medskip

\emph{Theorem C:
 For $v$ satisfying
$$a^*_{ij}\partial_{ij}v\ge 0, \quad \lambda I\le (a^*_{ij}(x))\le \Lambda I  $$
and
$$v\le 1\quad \mbox{ in } B_1 \quad \mbox{ and } \quad |\{v\le 1-\epsilon\}\cap B_1|\ge \mu |B_1|, $$
then $v\le 1-c(n,\lambda,\Lambda,\epsilon,\mu)$ over $B_{1/2}$
}
\medskip

By (\ref{delta-vij}) we see that
$$ \tilde a_{ij}(\lambda \cdot)\partial_{ij}(\Delta_{\hat e}^2v_{\lambda}+h_{\lambda})\ge 0, \quad \mbox{in }\quad B_1. $$
Applying Theorem C to
$\Delta_{\hat e}^2v_{\lambda}+h_{\lambda}$ we have
\begin{equation}\label{2e-3}
\Delta_{\hat e}^2v_{\lambda}+h_{\lambda}\le \alpha-C,\quad \mbox{in }\quad B_{1/2}.
\end{equation}
Since $\Delta_{\hat e}^2v_{\lambda}+h_{\lambda}$ is sub-harmonic,
$$\alpha=\sup_{\mathbb R^n}\Delta_e^2v+h=\lim_{\lambda\to \infty}\sup_{B_{1/2}}(\Delta_{\hat e}^2v_{\lambda}+h_{\lambda})<\alpha. $$
Then we get a contradiction, (\ref{2-e}) is established.

\medskip

\noindent{\bf Step four: The proof of the Liouville theorem for $n\ge 4$.}

\medskip

By a result of Li \cite{li-cpam}
there exists $\xi\in C^{2,\alpha}(\mathbb R^n)$ such that
$$\left\{\begin{array}{ll}
det(I+D^2\xi)=f_{v,p},\quad \mbox{ in }\quad \mathbb R^n, \\
I+D^2\xi>0,
\end{array}
\right.
$$
where $\xi$ is a periodic function with the same period as that of $f_{v,p}$.
Let $P(x)=\frac 12|x|^2$ and
$$w(x)=v(x)-P(x)-\xi(x). $$
Using $\xi(x+e)=\xi(x)$ for all $x\in \mathbb R^n$ and $e\in E$ , $\Delta_e^2P=1$ and (\ref{2-e}) we have
\begin{equation}\label{al-ca}
\Delta_e^2w+h\le 0, \quad \forall x\in \mathbb R^n, \quad \forall e\in E.
\end{equation}

From the equations for $v$ and $P+\xi$ we have,
$$det(D^2v)-det(D^2(P+\xi))=f_v-f_{v,p}, \quad \mbox{in}\quad \mathbb R^n, $$
which gives
$$
\hat a_{ij}\partial_{ij}w=f_v-f_{v,p}\quad \mbox{in }\quad \mathbb R^n
$$
where
$$\hat a_{ij}(x)=\int_0^1 cof_{ij}(tD^2v+(1-t)D^2(P+\xi))dt $$
is uniformly elliptic and divergence free. Using the Green's function corresponding to $-\partial_i(\hat a_{ij}(x)\partial_j)$ we
find $h_3$ to solve
$$\hat a_{ij}\partial_{ij}h_3=f_{v,p}-f_v,\quad \mbox{in}\quad \mathbb R^n $$
and
$$h_3(x)=O((1+|x|)^{-\beta_1}), \quad x\in \mathbb R^n. $$
Therefore we have
\begin{equation}\label{w-sol}
\hat a_{ij}\partial_{ij}(w+h_3)=0,\quad \mbox{in }\quad \mathbb R^n, \quad n\ge 3.
\end{equation}
Next we let
$$M_k=\sup_{B(0,k)}(w+h_3). $$
Then we claim that
\begin{equation}\label{double-M}
M_k\le 2M_{k/2}+C.
\end{equation}
for $C>0$ independent of $k$.

In order to prove (\ref{double-M}) we consider the equation for $\Delta_e^2w$ where $e\in E$.
In view of (\ref{p3e2}) we have
$$\tilde a_{ij}\partial_{ij}(\Delta_e^2w)\ge E_2,\quad \mbox{ in }\quad \mathbb R^n. $$
where $\tilde a_{ij}(x)=cof_{ij}(D^2v)$ is divergence free and $|E_2(x)|\le C(1+|x|)^{-2-\beta}$.

Just like the estimate of $h$ we can find $h_2$ that satisfies
$$\tilde a_{ij}\partial_{ij}(\Delta_e^2w+h_2)\ge 0 \quad \mbox{ in }\quad \mathbb R^n $$
and
$$|h_2(x)|\le C(1+|x|)^{-\beta_1},\quad \mbox{ in }\quad \mathbb R^n. $$

Let $x_0$ be where $\Delta_e^2 w+h_2$ attains its maximum on $B_k$, clearly $x_0\in \partial B_k$. Let $e=\frac{x_0}2+a$ be chosen so that $e\in E$ and $|a|$ is small:
$|a|\le C$.
Then for $k$ large
$$\Delta_e^2w(x-e)\le Ck^{-\beta_1}$$
which is
\begin{equation}\label{lconv}
w(x_0)\le 2w(x_0-e)-w(x_0-2e)+ O(k^{2-\beta_1}).
\end{equation}
Since $|x_0-2e|$ is in the neighborhood of $0$, we have $|w(x_0-2e)|\le C$
and
$$w(x_0)\le 2\max_{B(0,k/2)}w+C. $$
Hence (\ref{double-M}) is verified.

Set
$$g_k(y)=(w(ky)+h_3(ky))/M_k,\quad |y|\le 1. $$

It follows from (\ref{double-M}) , $w_1(0)=0$ and $h_3=O(1)$ that
$$\max_{B_1}g_k\to 1,\quad \max_{B_{1/2}}g_k\ge \frac 14, \quad g_k(0)\to 0, \mbox{ as } \quad k\to \infty. $$
Then we observe that the Harnack inequality holds for $1-g_k$ in $B_1$ because $1-g_k\ge 0$ in $B_1$ and
$$\hat a_{ij}(k\cdot)\partial_{ij}(1-g_k)=0 \quad \mbox{ in }\quad B_1. $$
The uniform ellipticity of $(\hat a_{ij})$ implies that
\begin{equation}\label{har-f}
\max_K(1-g_k)\le C(K)\min(1-g_k),\quad \forall K\subset\subset B_1.
\end{equation}
Consequently $g_k$ converges in $C^{\alpha}$ norm to $g$ in $B_1$ for some $\alpha>0$ depending only on the upper bound and the lower bound of the eigenvalues of $(\hat a_{ij})$. Since the constant in (\ref{har-f}) does not depend on $k$, we also have
$$\max_K(1-g)\le C(K)\min_K(1-g),\quad \forall K\subset\subset B_1. $$

Next we observe that (\ref{lconv}) also holds for all $e\in E$ with $\|e\|\le \frac 34k$, thus
$$g \mbox{ is concave in } B_1 $$
because the perturbation term disappears as $k\to \infty$.
Let $l$ be a linear function that touches $g$ from above around $0$ in $B_{1/2}$. Then $l-g$ is a convex function that takes its minimum at the origin. Since $l-g$ satisfies Harnack inequality, $l-g\equiv 0$ in $B_1$ (by the $C^{\alpha}$ convergence from $g_k$ to $g$, one can find $l_k\to l$ such that $l_k-g_k\ge 0$ in $B_{3/4}$ and $(l_k-g_k)(0)\to 0$. Thus applying Harnack inequality to $l_k-g_k$ we see that $l-g$ satisfies Harnack inequality too.

Since $\max_{B_1}g\ge \frac 34$ we see that $l=a\cdot x$ for $a\neq 0$, $|a|\ge 3/4$.  In other words
$$\frac{w(ky)}{M_k}-ay=o(1),\quad |y|\le \frac 34. $$
 Let $e\in E$ be a direction that $w(ke)\ge \frac 12 M_k k\|e\|$ for $k$ large. Subtract a linear function from $w$ ($w_1=w-$ the linear function) to make
$$w_1(0)=w_1(e)=0. $$
Clearly
\begin{equation}\label{w1-a}
w_1(ke)\ge \frac 14 M_k k\|e\|
\end{equation}
for $k$ large. Using
$$\Delta_e^2w_1(x)\le c_0|x|^{-\beta_1}, \quad |x|\ge 1 $$
we have
$$w_1(2e)\le 2w_1(e)-w_1(0)+c_0\|e\|^2\|e\|^{-\beta_1}=c_0 \|e\|^{2-\beta_1}. $$
Similarly
\begin{align*}
w_1(4e)&\le 2w_1(2e)-w_1(0)+c_0\|2e\|^2\|2e\|^{-\beta}\\
&\le 2c_0\|e\|^{2-\beta_1}+c_0(2\|e\|)^{2-\beta_1}.
\end{align*}

In general
\begin{align}\label{w1-b}
&w_1(2^{N+1}e)\\
\le &\|e\|^{2-\beta_1}c_0(2^{N}+2^{N-1+(2-\beta_1)}+2^{N-2+2(2-\beta_1)}+...+2^{N(2-\beta_1)}) \nonumber\\
= &\|e\|^{2-\beta}c_0 2^N\frac{1-2^{-(N+1)(\beta_1-1)}}{1-2^{1-\beta_1}}\nonumber
\end{align}
if $\beta_1>1$. Here we recall that $\beta_1=\min\{n-2,\beta\}$. Thus for $n\ge 4$ we obviously have
$$w_1(2^{N+1}e)\le C\|2^{N+1}e\|, $$
which leads to a contradiction (\ref{w1-a}). Therefore we have proved that $M_k$ is bounded if the Monge-Amp\`ere equation is defined in $\mathbb R^n$ for $n\ge 4$.
Since $w+h_3$ is bounded from above, it is easy to see that $w+h_3\equiv C$ by Harnack inequality.
 Theorem \ref{thm1} is established for $n\ge 4$.

 \medskip

 \noindent{\bf Step five: Proof of Theorem \ref{thm1} for $n=3$}

 \medskip

 In $\mathbb R^3$ we have $\beta_1=1$ which means (\ref{w1-b}) becomes
$$w_1(2^{N+1}e)\le CN(2^{N+1}\|e\|) $$
which leads to a contradiction if we assume $M_k\ge (\log k)^2$.
Therefore the same argument has proved that
$$w(x)+h_3(x)\le C(\log (2+|x|))^2, \quad x\in \mathbb R^3. $$

In order to finish the proof we need the following

\begin{lem}\label{har-4}
Let $u$ solve
$$a^*_{ij}\partial_{ij}u=0,\quad \mbox{in }\quad \mathbb R^n $$
where $\lambda I\le (a^*_{ij}(x))\le \Lambda I$ for all $x\in \mathbb R^n$ and
$$|u(x)|\le C(1+|x|)^{\delta}, \quad x\in \mathbb R^n. $$
There exists $\epsilon_0(\lambda,\Lambda)>0$ such that if $\delta\in (0,\epsilon_0)$
$u\equiv constant$.
\end{lem}

\noindent{\bf Proof of Lemma \ref{har-4}}:
For $R>1$ let
$$u_R(y)=u(Ry)/R^2, \quad |y|\le 1. $$
Then $u_R$ satisfies
$$a^*_{ij}(Ry)\partial_{ij}u_R=0,\quad \mbox{ in }\quad B_1 $$
and
$$|u_R(y)|\le CR^{\delta-2}, \quad \mbox{in }\quad B_1. $$
By Krylov-Safonov's \cite{safonov} estimate
\begin{equation}\label{safa-k}
\frac{|u_R(y)-u_R(0)|}{|y|^{\epsilon_0}}\le CR^{\delta}/R^2,\quad \forall y\in B_{1/2}
\end{equation}
where $\epsilon_0>0$ only depends on $\lambda$ and $\Lambda$.
Clearly (\ref{safa-k}) can be written as
$$|u(x)-u(0)|\le C|x|^{\epsilon_0}R^{\delta-\epsilon_0}, \quad \forall x\in B_{R/2}. $$
Fix any $x\in \mathbb R^n$, we let $R\to \infty$, then $u(x)=u(0)$ for $\delta<\epsilon_0$.
Lemma \ref{har-4} is established. $\Box$.

\medskip

Applying Lemma \ref{har-4} to $w+h_3$ we see that $w+h_3\equiv c$ in $\mathbb R^3$. Theorem \ref{thm1} is also proved for the equation defined in $\mathbb R^3$. $\Box$


\begin{thebibliography}{99}

\bibitem{alek1} A. D. Aleksandrov, Dirichlet's problem for the equation $Det \|z_{ij}\|=\phi$ I, Vestnik Leningrad. Univ. Ser. Mat. Meh. Astr. 13 (1958), 5--24.
\bibitem{b1} I. J. Bakelman, Generalized solutions of Monge-Amp\`ere equations (Russion). Dokl. Akad. Nauk SSSR (N. S. ) 114 (1957), 1143--1145.
\bibitem{bao} J. Bao, H. Li, On the exterior Dirichlet problem for the Monge-Amp\`ere equation in dimension two,  Nonlinear Anal. 75 (2012), no. 18, 6448–-6455.
\bibitem{zhang-mp} J. Bao, H. Li, L. Zhang, Monge-Amp\'ere equations on exterior domains, Cal. Var. PDE. in press, arxiv: 1304:2415

\bibitem{caf-ann1} L. A. Caffarelli, A localization property of viscosity solutions to the Monge–Amp\`ere equation and
their strict convexity, Ann. of Math. 13 (1990) 129-–134.

\bibitem{caf-ann2} L. A. Caffarelli, Interior $W^{2,p}$ estimates for solutions of the Monge–Amp\`ere equation, Ann. of
Math. 131 (1990) 135-–150.
\bibitem{caf-cab} L. A. Caffarelli, X. Cabre, Fully nonlinear elliptic equations. American Mathematical Society Colloquium Publications, 43. American Mathematical Society, Providence, RI, 1995. vi+104 pp. ISBN: 0-8218-0437-5.

\bibitem{caf-li-1} L. A. Caffarelli, Y. Y. Li, An extension to a theorem of J\"orgens, Calabi, and Pogorelov. Comm. Pure Appl. Math.
56 (2003), no. 5, 549-–583.
\bibitem{caf-li-2} L. A. Caffarelli, Y. Y. Li, A Liouville theorem for solutions of the Monge-Amp\`ere equation with periodic data. Ann. Inst. H. Poincaré Anal. Non Linéaire 21 (2004), no. 1, 97-–120.
\bibitem{gute}  L. A. Caffarelli, C. E. Gutiérrez. Properties of the solutions of the linearized Monge-Ampère equation. Amer. J. Math. 119 (1997), no. 2, 423–465.
\bibitem{cns} L. A. Caffarelli, L. Nirenberg, J. Spruck, The Dirichlet problem for nonlinear second-order elliptic equations. I. Monge-Amp\`ere equation. Comm. Pure Appl. Math. 37 (1984), no. 3, 369–-402.
\bibitem{caf-via} L. A. Caffarelli, J. Viaclovsky, On the regularity of solutions to Monge-Ampère equations on Hessian manifolds.
Comm. Partial Differential Equations 26 (2001), no. 11-12, 2339–2351.

\bibitem{calabi} E. Calabi, Improper affine hyperspheres of convex type and a generalization of a theorem by K. J\"orgens.
Michigan Math. J. 5 (1958) 105-–126.
\bibitem{cheng-yau1} S. Y. Cheng, S. T. Yau, On the regularity of the Monge-Amp\`ere equation $\det(\partial^2u/\partial x_i\partial x_j)=F(x,u)$.
Comm. Pure Appl. Math. 30 (977), no. 1, 41--68.
\bibitem{cheng-yau2} S. Y. Cheng, S. T. Yau, Complete affine hypersurfaces. I. The completeness of affine metrics.
Comm. Pure Appl. Math. 39 (1986), no. 6, 839-–866.
\bibitem{chou-wang} K. S. Chou, X. Wang, Entire solutions of the Monge-Amp\`ere equation.
Comm. Pure Appl. Math. 49 (1996), no. 5, 529–-539.

\bibitem{delanoe} P. Delano\"e, Partial decay on simple manifolds. Ann. Global Anal. Geom. 10 (1992), no. 1, 3--61.
\bibitem{figalli}  G. De Philippis, A. Figalli, $W^{2,1}$ regularity for solutions of the Monge-Ampère equation. Invent. Math. 192 (2013), no. 1, 55–-69.

\bibitem{savin2}  G. De Philippis, A. Figalli, O. Savin,  A note on interior W2,1+e estimates for the Monge-Ampère equation. Math. Ann. 357 (2013), no. 1, 11–-22.

\bibitem{FMM1} L, Ferrer, ; A. Mart\'inez, ; F. Mil\'an,  The space of parabolic affine spheres with fixed compact boundary. Monatsh. Math. 130 (2000), no. 1, 19–-27.
\bibitem{FMM2} L, Ferrer, ; A. Mart\'inez, ; F. Mil\'an, An extension of a theorem by K. J\"orgens and a maximum principle at infinity for parabolic affine spheres. Math. Z. 230 (1999), no. 3, 471–-486.

\bibitem{gilbarg} D. Gilbarg; Neil S. Trudinger, Elliptic partial differential equations of second order. Reprint of the 1998 edition. Classics in Mathematics. Springer-Verlag, Berlin, 2001. xiv+517 pp. ISBN: 3-540-41160-7.

\bibitem{huang} Q. Huang, Sharp regularity results on second derivatives of solutions to the Monge Ampère equation with VMO type data. Comm. Pure Appl. Math. 62 (2009), no. 5, 677–-705.
\bibitem{ivo-1} Ivochkina, N.M.: Construction of a priori estimates for convex solutions of the Monge Ampère equation by the integral method (Russian). Ukrain. Mat. Z. 30, 45–-53 (1978)
\bibitem{ivo-2} Ivochkina, N.M.: A priori estimate of $\|u\|_{C^{2,\alpha}(\Omega)}$ of convex solutions of the Dirichlet problem for the Monge Ampère equation. Zap. Nauchn. Sem. Leningrad. Otdel. Mat. Inst. Steklov. 96, 69–79 (1980).
\bibitem{ivo-3} Ivochkina, N.M.: Classical solvability of the Dirichlet problem for the Monge–Ampère equation (Russian). Zap. Nauchn. Sem. Leningrad. Otdel. Mat. Inst. Steklov (LOMI) 131, 72–79 (1983).

\bibitem{jian-wang} H. Jian, X. Wang, Continuity estimates for the Monge-Ampère equation. SIAM J. Math. Anal. 39 (2007), no. 2, 608–-626.

\bibitem{jorgens} K, J\"orgens, \"Uber die L\"osungen der Differentialgleichung $rt-s^2=1$. Math. Ann. 127 (1954) 130--134.

\bibitem{jost-xin} J. Jost and Y.L. Xin, Some aspects of the global geometry of entire space-like
submanifolds. Results Math. 40 (2001), 233--245.

\bibitem{k2} N. V. Krylov, Boundedly inhomogeneous elliptic and parabolic equations in a domain. (Russian) Izv. Akad. Nauk SSSR Ser. Mat. 47 (1983), no. 1, 75–-108.

\bibitem{safonov} N. V. Krylov, M. V. Safonov, An estimate for the probability of a diffusion process hitting a set of positive measure, Doklady Akademii Nauk SSSR 245 (1): 18–20, ISSN 0002-3264

\bibitem{li-cpam} Y. Y. Li, Some existence results of fully nonlinear elliptic equations of Monge–Ampère type, Comm. Pure Appl. Math. 43 (1990) 233–271.
\bibitem{littman} W. Littman, G. Stampacchia, H. F. Weinberger,
Regular points for elliptic equations with discontinuous coefficients.
Ann. Scuola Norm. Sup. Pisa (3) 17 1963 43–-77.

\bibitem{n2} L. Nirenberg, The Weyl and Minkowski problems in differential geometry in the large. Comm. Pure Appl. Math. 6 (1953), 337-394.

\bibitem{p5} A. V. Pogorelov, The regularity of the generalized solutions of the equation $det(\partial^2u/\partial x_i\partial x_j)=\phi>0$, (Russian) Dokl. Akad. Nauk SSSR 200 (1971), 534–-537.

\bibitem{savin-1} O. Savin, A Liouville theorem for solutions to the linearized Monge-Ampere equation. Discrete Contin. Dyn. Syst. 28 (2010), no. 3, 865–873.
\bibitem{tw7} N. Trudinger, X. Wang, Boundary regularity for the Monge-Ampère and affine maximal surface equations.
Ann. of Math. (2) 167 (2008), no. 3, 993–-1028.

\bibitem{tw1} N. Trudinger, X. Wang, The Monge-Ampère equation and its geometric applications. Handbook of geometric analysis. No. 1, 467–-524, Adv. Lect. Math. (ALM), 7, Int. Press, Somerville, MA, 2008.

\bibitem{urbas} J. Urbas, Regularity of generalized solutions of Monge-Amp\`ere equations.
Math. Z. 197 (1988), no. 3, 365–-393.

\end{thebibliography}
\end{document}